\documentclass[12pt]{amsart}
\usepackage{amssymb, amsmath, url, graphicx}

\newtheorem{theorem}{Theorem}

\theoremstyle{definition}
\newtheorem{remark}[theorem]{Remark}

\newcommand{\Z}{\mathbb{Z}}
\newcommand{\dd}{\displaystyle }

\title{The subgroup commutativity degree of finite $P$-groups}

\author{Marius T\u arn\u auceanu}
\address{Marius T\u arn\u auceanu \\
Faculty of  Mathematics \\
''Al.I. Cuza'' University \\
Ia\c si \\
Romania}
\email{tarnauc@uaic.ro}

\date{December 1, 2013}

\begin{document}

\begin{abstract}
    The subgroup commutativity degree of a group $G$ has been defined in \cite{6} as the probability
    that two subgroups of $G$ commute, or equivalently that the product of two subgroups is again a subgroup.
    Problem 4.3 of \cite{6} asks whether there exist families of groups other than dihedral, quasi-dihedral or generalized quaternion (all of $2$-power cardinality), whose subgroup commutativity degree tends to $0$ as the size of the group tends to infinity. An affirmative answer to this question has been provided by Aivazidis \cite{1,2} for the family of projective special linear groups over fields of even characteristic and for the family of the simple Suzuki groups. In this short note we indicate another family of groups with this property, namely the finite $P$-groups.
\end{abstract}


\maketitle

\section{Introduction}

\textit{The subgroup commutativity degree} (also called \textit{the subgroup permutability degree}) of a finite group $G$ is defined as the probability that two subgroups of $G$ commute, or equivalently that the product of two subgroups is again a subgroup
\begin{equation}
sd(G)=\frac{1}{|L(G)|^2}\left|\{(H,K){\in} L(G)^2\mid HK=KH\}\right|,\nonumber
\end{equation}where $L(G)$ is the subgroup lattice of $G$. Many problems related to this concept have been formulated in \cite{6,8}. We recall here only the problem of finding some natural families of groups $G_n$, $n\in\mathbb{N}$, whose subgroup commutativity degree vanishes asymptotically, i.e.
\begin{equation}
\lim_{n\to\infty}sd(G_n)=0. \tag{1}
\end{equation}It is known that the dihedral groups $D_{2^n}$, the quasi-dihedral groups $S_{2^n}$, the generalized quaternion groups $Q_{2^n}$, the projective special li\-ne\-ar groups ${\rm PSL}_2(2^n)$ and the simple Suzuki groups ${\rm Sz}(2^{2n+1})$ satisfy (1). Our main result shows that the finite $P$-groups also satisfy this property.

\begin{theorem}\label{main}
    The subgroup commutativity degree of\, $G_{n,p}$ vanishes asym\-pto\-ti\-ca\-lly, i.e.
    \begin{equation}
    \lim_{n\to\infty}sd(G_{n,p})=0. \nonumber
    \end{equation}
\end{theorem}

Most of our notation is standard and will usually not be repeated
here. Basic definitions and results on groups can be found in \cite{4}.
For subgroup lattice concepts we refer the reader to \cite{3}.

\section{Preliminaries}

In the following we will recall the notion of $P$-group, according
to \cite{3}, and we will indicate the subgroup structure of these groups.

Let $p$ be a prime, $n\geq 2$ be a cardinal number
and $G$ be a group. We say that $G$ belongs to the class $P(n,p)$
if it is either elementary abelian of order $p^n$, or a semidirect
product of an elementary abelian normal subgroup $H$ of order $p^{n-1}$
by a group of prime order $q\neq p$ which induces a nontrivial power
automorphism on $H$. The group $G$ is called a \textit{P-group} if
$G\in P(n,p)$ for some prime $p$ and some cardinal number $n\geq 2$.
It is well-known that the class $P(n,2)$ consists only of the elementary
abelian group of order $2^n$. Also, for $p>2$ the class $P(n,p)$ contains
the elementary abelian group of order $p^n$ and, for every prime
divisor $q$ of $p-1$, exactly one nonabelian $P$-group with
elements of order $q$. Moreover, the order of this group is
$p^{n-1}q$ if $n$ is finite. The most important property of the
groups in a class $P(n,p)$ is that they are all lattice-isomorphic
(see Theorem 2.2.3 of \cite{3}). This played an essential role in
\cite{6} to produce examples of finite lattice-isomorphic groups with
different subgroup commutativity degrees.

Since the subgroup commutativity degree concept is defined only for finite
groups and it is trivial in the abelian case, we will focus only on finite
nonabelian $P$-groups. So, let us suppose that $p>2$ and $n\in\mathbb{N}$ are
fixed, and take a divisor $q$ of $p-1$. The nonabelian group of
order $p^{n-1}q$ in the class $P(n,p)$ will be denoted by
$G_{n,p}$\,. By Remarks 2.2.1 of \cite{3}, it is of type
\begin{equation}
G_{n,p}=H\langle x\rangle,\nonumber
\end{equation}where $H\cong\mathbb{Z}_p^{n-1}$ (i.e. the direct product of $n-1$ copies
of $\mathbb{Z}_p$), $o(x)=q$ and there exists an integer $r$ such that
$x^{-1}hx=h^r$, for all $h\in H$.

In order to describe the subgroups of $G_{n,p}$\,,
we need some information about the subgroups of a finite
ele\-men\-tary abelian $p$-group. First of all, we recall the following
well-known theorem (see e.g. \cite{5,7}).

\begin{theorem}\label{Sd}
    The number of all subgroups of order $p^k$ of the finite elementary
    abelian $p$-group $\mathbb{Z}_p^n$ is $1$ if $k=0$ or $k=n$, and
    \begin{equation}
    a_{n,p}(k)=\displaystyle\sum_{1\le i_1<i_2<...<i_k\le n}p^{i_1+i_2+...+i_k-\frac{k(k+1)}2}\nonumber
    \end{equation}if $1\le k\le n-1$. In particular, the total number of subgroups of $\mathbb{Z}_p^n$ is
    \begin{equation}
    a_{n,p}=2+\displaystyle\sum_{k=1}^{n-1}a_{n,p}(k)=2+\displaystyle\sum_{k=1}^{n-1} \displaystyle\sum_{1\le i_1<i_2<...<i_k\le n}p^{i_1+i_2+...+i_k-\frac{k(k+1)}2}.\nonumber
    \end{equation}
\end{theorem}Note that an alternative way of writing the numbers $a_{n,p}(k)$, $k=0,1,...,n$, is
\begin{equation}
a_{n,p}(k)=\frac{(p^n-1)\cdots (p-1)}{(p^k-1)\cdots (p-1)(p^{n-k}-1)\cdots (p-1)}\,.\nonumber
\end{equation}

Some properties of the above numbers will be very useful in
determining $sd(G_{n,p})$.

\begin{remark}
    The numbers $a_{n,p}(k)$ and $a_{n,p}$ in Theorem 2 satisfy the following recurrence relations:
\begin{itemize}
        \item[\rm 1)]
        $a_{n,p}(k)=a_{n-1,p}(k)+p^{n-k}a_{n-1,p}(k-1)$, for all $k=1,2,...,n-1$,
        \item[\rm 2)]
        $a_{n,p}=2a_{n-1,p}+(p^{n-1}-1)a_{n-2,p}$\,.
\end{itemize}We also remark that $a_{n,p}$ can be written as $a_{n,p}=f_n(p)$,
where the polynomial $f_n\in{\rm\mathbb{Z}[X]}$ is of degree
$\left[\frac{n^2}4\right]$. The dominant coefficient $x_n$ of this
polynomial is 1 if $n$ is even and 2 if $n$ is odd. Moreover, by
using a computer algebra program, from 2) we can easily obtain the
first terms of the integer sequence $(a_{n,p})_{n\in\mathbb{N}^*}$. For example, we have:
\begin{equation}
\begin{array}{l}
a_{1,p}=2,\\
a_{2,p}=p+3,\\
a_{3,p}=2p^2+2p+4,\\
a_{4,p}=p^4+3p^3+4p^2+3p+5,\\
a_{5,p}=2p^6+2p^5+6p^4+6p^3+6p^2+4p+6,\\
\vdots\\
\end{array}\nonumber
\end{equation}and so on.
\end{remark}

A subgroup of $G_{n,p}$ is either cyclic if it is included in $H$,
or a semidirect product of the same type as $G_{n,p}$ if it
possesses some elements of order $q$. So, we can give an
enumerative description of these subgroups. They are:
\begin{itemize}
        \item[\rm -]
        one of order 1, namely the trivial subgroup $H^1_1$,
        \item[\rm -]
        $a_{n-1,p}(1)$ of order $p$, say $H^p_i$, $i=1, ..., a_{n-1,p}(1)$,
        \item[\rm -]
        $a_{n-1,p}(2)$ of order $p^2$, say $H^{p^2}_i$, $i=1, ..., a_{n-1,p}(2)$,\\
        \vdots
        \item[\rm -]
        $a_{n-1,p}(n-2)$ of order $p^{n-2}$, say $H^{p^{n-2}}_i$, $i=1, ..., a_{n-1,p}(n{-}2)$,
        \item[\rm -]
        one of order $p^{n-1}$, namely $H^{p^{n-1}}_1=H$,
        \item[\rm -]
        $p^{n-1}$ of order $q$, say $H^q_i$, $i=1, ..., p^{n-1}$,
        \item[\rm -]
        $a_{n-1,p}(1)p^{n-2}$ of order $pq$, say $H^{pq}_i$, $i=1, ..., a_{n-1,p}(1)p^{n-2}$,
        \item[\rm -]
        $a_{n-1,p}(2)p^{n-3}$ of order $p^2q$, say $H^{p^2q}_i$, $i=1, ..., a_{n-1,p}(2)p^{n-3}$,\\
        \vdots
        \item[\rm -]
        $a_{n-1,p}(n{-}2)p$ of order $p^{n-2}q$, say $H^{p^{n-2}q}_i$, $i=1, ..., a_{n-1,p}(n{-}2)p$,
        \item[\rm -]
        one of order $p^{n-1}q$, namely $H^{p^{n-1}q}_1=G_{n,p}$\,.
\end{itemize}

We observe that
\begin{equation}
|L(G_{n,p})|=a_{n,p}\nonumber
\end{equation}since $G_{n,p}$ and $\mathbb{Z}_p^n$ are lattice-isomorphic. On the other hand, by Lemma 2.2.2 of \cite{3} we infer that the normal subgroups of $G_{n,p}$ are $G_{n,p}$ itself and all subgroups contained in $H$. Therefore
\begin{equation}
|N(G_{n,p})|=1+|L(H)|=1+|L(\Z^{n-1}_p)|=1+a_{n-1,p}\,,\nonumber
\end{equation}where $N(G_{n,p})$ denotes the normal subgroup lattice of $G_{n,p}$\,.

We are now able to prove our main result.

\section{Proof of Theorem 1}

First of all, we will prove the following inequality
\begin{equation}
sd(G_{n,p})\leq \frac{a_{n-1,p}}{a_{n,p}}\left(2+\frac{1}{a_{n,p}}\right).\tag{2}
\end{equation}

For every subgroup $K$ of $G_{n,p}$, let us denote by $C(K)$
the set of all subgroups of $G_{n,p}$ which commute with $K$. Then
\begin{equation}
sd(G_{n,p})=\frac{1}{|L(G_{n,p})|^2}\,\dd\sum_{K\in L(G_{n,p})}|C(K)|=\frac{1}{a_{n,p}^2}\,\dd\sum_{K\in L(G_{n,p})}|C(K)|\,.\tag{3}
\end{equation}Moreover, we have
\begin{equation}
|C(H^{p^k}_i)|= a_{n,p}\,, \forall\, k=0,1,...,n-1 \mbox{ and } \forall\, i=1,2,...,a_{n-1,p}(k),\nonumber
\end{equation}because all $p$-subgroups of $G_{n,p}$
are normal. Then (3) becomes
\begin{equation}
sd(G_{n,p})=\frac{1}{a_{n,p}^2}\left(a_{n-1,p}\,a_{n,p}+\dd\sum_{k=0}^{n-1}\,\sum_{i=1}^{a_{n-1,p}(k)}|C(H^{p^kq}_i)|\right).\tag{4}
\end{equation}

Assume that $k\in\{0,1, ..., n-1\}$ and $i\in\{1,2, ...,
a_{n-1,p}(k)\}$ are fixed, and take a subgroup $S\in
C(H^{p^kq}_i)$. Then either $S\in N(G_{n,p})$ or
$q\hspace{1mm}|\hspace{1mm}|S|$. In the second case, by the equality
\begin{equation}
|SH^{p^kq}_i|=\dd\frac{|S||H^{p^kq}_i|}{|S\cap H^{p^kq}_i|}\nonumber
\end{equation}it follows that $q$ must divide $|S\cap H^{p^kq}_i|$ and so there is
a subgroup $H^q$ of order $q$ of $G_{n,p}$ contained both in $S$
and $H^{p^kq}_i$. Thus
\begin{equation}
C(H^{p^kq}_i)=N(G_{n,p})\cup \left(\dd\bigcup_{H^q\in Q} \{S\in
L(G_{n,p})\,|\,H^q \subseteq S\}\right)\nonumber
\end{equation}
\begin{equation}
\hspace{9mm}=N(G_{n,p})\cup\left(\dd\bigcup_{H^q\in Q} \{H^qT \,|\,T\in L(H)\}\right),\nonumber
\end{equation}where $Q$ denotes the set of all subgroups of order $q$ in
$H^{p^kq}_i$. Since $|Q|=p^k$, one obtains
\begin{equation}
|C(H^{p^kq}_i)|\leq |N(G_{n,p})|+|L(H)|p^k=1+a_{n-1,p}\,(1+p^k).\nonumber
\end{equation}This implies that
\begin{equation}
\dd\sum_{k=0}^{n-1}\,\sum_{i=1}^{a_{n-1,p}(k)}|C(H^{p^kq}_i)|\leq
a_{n-1,p}\left(1+\dd\sum_{k=0}^{n-1}\,\sum_{i=1}^{a_{n-1,p}(k)}(1+p^k)\right)\nonumber
\end{equation}
\begin{equation}
=a_{n-1,p}\left(1+a_{n-1,p}+\dd\sum_{k=0}^{n-1}a_{n-1,p}(k)p^k\right)=
a_{n-1,p}\left(1+a_{n,p}\right),\nonumber
\end{equation}where the last equality has been obtained from the recurrence
relation 1). Then (4) shows that
\begin{equation}
sd(G_{n,p})\leq \frac{a_{n-1,p}}{a_{n,p}}\left(2+\frac{1}{a_{n,p}}\right),\nonumber
\end{equation}as desired.

Since $a_{n,p}$ can be written as a
polynomial in $p$ of degree $\left[\frac{n^2}4\right]$ and
dominant coefficient $x_n \in \{1,2\}$, we have
\begin{equation}
\dd\lim_{n\to\infty}\frac{a_{n-1,p}}{a_{n,p}}=\dd\lim_{n\to\infty}\frac{x_{n-1}}{x_n}\hspace{1mm}p^{\hspace{0,5mm}\left[\frac{(n-1)^2}4\right]-\left[\frac{n^2}4\right]}=\dd\lim_{n\to\infty}\frac{x_{n-1}}{x_n}\hspace{1mm}p^{-\left[\frac{n}2\right]}=0,\nonumber
\end{equation}which together with (2) lead to
\begin{equation}
\dd\lim_{n\to\infty}sd(G_{n,p})=0.\nonumber
\end{equation}This completes the proof.
\hfill\rule{1,5mm}{1,5mm}

\end{document}